\documentclass[11pt,reqno]{amsart}

\usepackage{amsmath, amsfonts, amsthm, amssymb}

\newtheorem{Theorem}{Theorem}[section]
\newtheorem{Lemma}{Lemma}[section]

\theoremstyle{definition}

\theoremstyle{remark}
\newtheorem{Remark}{Remark}[section]

\usepackage{amsmath}

\renewcommand{\u}{{\bf u}}

\newcommand{\R}{{\mathbb R}}
\newcommand{\Dv}{{\rm div}}

\def\f{\frac}
\renewcommand{\O}{\Omega}

\def\D{\Delta }

\def\hf1{^\f{1}{1-\xi^2}}

\def\be{\begin{equation}}
\def\en{\end{equation}}
\def\bs{\begin{split}}
\def\es{\end{split}}

\author{Cheng Yu}
\address{Department of Mathematics,  The University of Texas,
                           Austin, Texas 78712.}
\email{yucheng@math.utexas.edu}

\title%[Global weak solutions to N-S-V equations]
[]
{ A new proof to the energy conservation for the Navier-Stokes equations}

\keywords{Energy conservation, Navier-Stokes equations, weak solution.}
\subjclass[2000]{}

\date{\today}

\begin{document}
\begin{abstract}
In this paper we give a new proof to the energy conservation for the weak solutions of the incompressible Navier-Stokes equations. This result was first proved by Shinbrot.
The new proof relies on a lemma introduced by Lions.
\end{abstract}

\maketitle
\section{Introduction}
We are interested in studying the energy conservation for the weak solutions of Navier-Stokes equations
\begin{equation}
\label{NS equations}
\begin{split}
& \u_t+\u\cdot\nabla\u+\nabla P-\mu\D\u=0,
\\&\Dv\u=0,
\end{split}
\end{equation}
with the initial data
\begin{equation}
\label{initial data}
\u(0,x)=\u_0
\end{equation}
for $(t,x)\in \R^{+}\times\O,$ where $\O=\mathbb{T}^d$ is a periodic domain in $\R^d.$
\vskip0.3cm

The existence of weak solution was proved by Leray \cite{Le} and Hopf \cite{H}. The notion of weak solution has been introduced  in \cite{Le}. As usual, a weak solution $\u$ satisfies the energy inequality
\begin{equation}
\label{energy inequality}
\int_{\O}|\u(t,x)|^2\,dx+2\mu\int_0^T\int_{\O}|\nabla\u|^2\,dx\,dt\leq \int_{\O}|\u_0|^2\,dx,
\end{equation}
for any $t\in(0,T).$
It is a natural question to ask when a weak solution satisfies the stronger version of \eqref{energy inequality}, that is,
\begin{equation}
\label{energy conservation}
\int_{\O}|\u(t,x)|^2\,dx+2\mu\int_0^T\int_{\O}|\nabla\u|^2\,dx\,dt=\int_{\O}|\u_0|^2\,dx.
\end{equation}
\vskip0.3cm

As we all known, any classical solution of the Navier-Stokes equations satisfies the energy equality \eqref{energy conservation}. However, the existence of global classical solution remains open. Thus, an interesting question is how badly behaved $\u$ can keep the energy conservation.  In his pioneering work \cite{Serrin}, Serrin has proved $\u$ satisfies \eqref{energy conservation} if $\u \in L^p(0,T;L^q(\O))$, where
\begin{equation}
\label{Serrion condition}
\frac{2}{p}+\frac{d}{q}\leq 1,
\end{equation}
where $d$ is the dimension of space.
In \cite{Shinbrot}, Shinbrot has shown the same conclusion if $\u \in L^p(0,T;L^q(\O))$, where
\begin{equation}
\label{Shinbrot condition}
\frac{1}{p}+\frac{1}{q}\leq \frac{1}{2},\; q\geq 4.
\end{equation}
Note that, it is hard to say which condition is weaker between \eqref{Serrion condition} and \eqref{Shinbrot condition}. However, an interesting point about the condition \eqref{Shinbrot condition} is that they do not depend on the dimension $d$.  We have to mention that a similar result to the Euler equations, which was proved by E-Constantin-Titi \cite{CET}.  It was the answer to the first part of Onsager's conjecture \cite{O}.

\vskip0.3cm

The goal of this paper is to give a new proof to Shinbrot's remarkable result in \cite{Shinbrot}.
The following is the main result of this paper:
\begin{Theorem}
\label{main result}
Let $\u\in L^{\infty}(0,T;L^2(\O))\cap L^2(0,T;H^1(\O))$ be a weak solution of the incompressible Navier-Stokes equations, that is,
\begin{equation}
\label{weak formulation}
\begin{split}
-\int_0^T\int_{\O}\u\varphi_t\,dx\,dt&-\int_{\O}\u_0\varphi(0,x)\,dx-\int_0^T\int_{\O}\nabla\varphi\u\otimes\u\,dx\,dt
\\&+\mu\int_0^T\int_{\O}\nabla\u\nabla\varphi\,dx\,dt=0
\end{split}
\end{equation}
for any smooth test function $\varphi\in C^{\infty}(\R^{+}\times\O)$ with compact support, and $\Dv\varphi=0.$ In additional, if $\u \in L^r(0,T;L^s(\O))$ for any $\frac{1}{r}+\frac{1}{s}\leq \frac{1}{2}, \; s\geq 4$, then
\begin{equation*}
%\label{energy conservation}
\int_{\O}|\u(t,x)|^2\,dx+2\mu\int_0^T\int_{\O}|\nabla\u|^2\,dx\,dt=\int_{\O}|\u_0|^2\,dx
\end{equation*}
for any $t\in [0,T].$
\end{Theorem}

\begin{Remark}
The proof of Theorem \ref{main result} was motivated
by the work of Vasseur-Yu \cite{VY}, where they have shown the first existence result of weak solutions to the degenerate compressible Navier-Stokes equations in dimension 3. The same conclusion for the compressible version is established in \cite{Yu}.
\end{Remark}

\begin{Remark} The global existence of weak solution to $2d$ Euler equations was proved in \cite{L}, in particular, see Theorem 4.1 of book \cite{L}. For a weak solution to Euler equation in this sense, adopting the same argument, we can conclude the energy conserve  for any weak solution $\u\in C([0,\infty);W^{1,r}(\O))$, where $r\geq \frac{3}{2}.$ It was mentioned on page 132 in \cite{L}.
\end{Remark}

\vskip0.3cm
\section{Proof}

The goal of this section is to prove our main result. To this end, we need to introduce a crucial lemma.
The key lemma is as follows which was proved by Lions in \cite{L}.
 \begin{Lemma}
 \label{Lions's lemma}
 Let $f\in W^{1,p}(\R^d),\,g\in L^{q}(\R^d)$ with $1\leq p,q\leq \infty$, and $\frac{1}{p}+\frac{1}{q}\leq 1$. Then, we have
 $$\|\Dv(fg)*\eta_{\varepsilon}-\Dv(f(g*\eta_{\varepsilon}))\|_{L^{r}(\R^d)}\leq C\|f\|_{W^{1,p}(\R^d)}\|g\|_{L^{q}(\R^d)}$$
 for some $C\geq 0$ independent of $\varepsilon$, $f$ and $g$, $r$ is determined by $\frac{1}{r}=\frac{1}{p}+\frac{1}{q}.$ In addition,
  $$\Dv(fg)*\eta_{\varepsilon}-\Dv(f(g*\eta_{\varepsilon}))\to0\;\;\text{ in }\,L^{r}(\R^d)$$
 as $\varepsilon \to 0$ if $r<\infty.$ Here  $\varepsilon>0$ is a small enough number, $\eta\in C_0^{\infty}(\O)$ be a standard mollifier supported in $B(0,1).$
 \end{Lemma}
The weak solution $\u$  is uniformly bounded in $L^{\infty}(0,T;L^2(\O))\cap L^2(0,T;H^1(\O))$. Thus, it is possible to make use of Lemma \ref{Lions's lemma} to handle convective term $\Dv(\u\otimes\u).$ With Lemma \ref{Lions's lemma} in hand, we are ready to prove our main result.

\vskip0.3cm
We define a new function $\Phi=\overline{\u}$, where $\overline{f(t,x)}=f*\eta_{\varepsilon}(x)$, $\varepsilon>0$ is a small enough number, $\eta\in C_0^{\infty}(\O)$ be a standard mollifier supported in $B(0,1).$ Note that, we have
\begin{equation}
\label{incompressibility}
\Dv\overline{\u}=0
\end{equation}
Using $\overline{\Phi}$ to test Navier-Stokes equations \eqref{NS equations}, ones obtain
\begin{equation*}
\int_{\O}\overline{\Phi}\left(\u_t+\Dv(\u\otimes\u)+\nabla P-\mu\Delta\u\right)\,dx=0,
\end{equation*}
which in turn gives us
\begin{equation*}
\int_{\O}\overline{\u}\overline{\left(\u_t+\Dv(\u\otimes\u)+\nabla P-\mu\Delta\u\right)}\,dx=0.
\end{equation*}
This yields
\begin{equation*}
\frac{1}{2}\frac{d}{dt}\int_{\O}|\overline{\u}|^2\,dx+\mu\int_{\O}|\nabla\overline{\u}|^2\,dx=\int_{\O}\overline{\Dv(\u\otimes\u)}\overline{\u}\,dx,
\end{equation*}
and hence
\begin{equation}
\label{energy equality before limit}
\begin{split}&
\int_{\O}|\overline{\u}|^2\,dx-
\int_{\O}|\overline{\u}_0|^2\,dx+2\mu\int_0^T\int_{\O}|\nabla\overline{\u}|^2\,dx\,dt
\\&\quad\quad\quad\quad\quad\quad\quad\quad=2\int_0^T\int_{\O}\overline{\Dv(\u\otimes\u)}\overline{\u}\,dx\,dt.
\end{split}
\end{equation}
Next we rewrite \begin{equation}
\label{convecation term}
\begin{split}
&\overline{\Dv(\u\otimes\u)}=\left(\overline{\Dv(\u\otimes\u)}-\Dv(\u\otimes\overline{\u})\right)
\\&+\left[\Dv(\u\otimes\overline{\u})
-\Dv(\overline{\u}\otimes\overline{\u})\right]
+\Dv(\overline{\u}\otimes\overline{\u})
\\&=R_1+R_2+R_3.
\end{split}
\end{equation}
Thus, the right-hand side of \eqref{energy equality before limit} is given by
$$\int_0^T\int_{\O}\left(R_1+R_2+R_3\right)\overline{\u}\,dx\,dt.$$
By means of \eqref{incompressibility}, we have
\begin{equation}
\label{R3 term}
\int_0^T\int_{\O}R_3\overline{\u}\,dx\,dt=0.
\end{equation}

\vskip0.3cm

Now we first assume that $\u\in L^p(0,T;L^q(\O)),$ where $p,q\geq 4$. This restriction will be improved at the very end.
We can control the term related to $R_2$ in the following way
\begin{equation}
\label{R2 term}
\begin{split}
&\left|\int_0^T\int_{\O}R_2\overline{\u}\,dx\,dt\right|
\\&=
\left|\int_0^T\int_{\O}(\u\otimes\overline{\u}-\overline{\u}\otimes\overline{\u})\nabla\overline{\u}\,dx\,dt\right|
\\&\leq  \int_0^T\int_{\O}|\u-\overline{\u}||\overline{\u}||\nabla\overline{\u}|\,dx\,dt
\\&\leq C\|\u-\overline{\u}\|_{L^p(0,T;L^q(\O))}\|\overline{\u}\|_{L^p(0,T;L^q(\O))}\|\nabla\overline{\u}\|_{L^2(0,T;L^2(\O))}
\\&\to0
\end{split}
\end{equation}
as $\varepsilon$ goes to zero, where $p,q\geq 4.$\\

Meanwhile, thanks to Lemma \ref{Lions's lemma}, we find
\begin{equation}
\label{convergence of R1}
\|R_1\|_{L^{\frac{2p}{p+2}}(0,T;L^{\frac{2q}{q+2}}(\O))} \leq C\|\overline{\u}\|_{L^p(0,T;L^q(\O))}\|\nabla\u\|_{L^2(0,T;L^2(\O))},
\end{equation}
and it converges to zero in  $L^{\frac{2p}{p+2}}(0,T;L^{\frac{2q}{q+2}}(\O))$ as $\varepsilon$ tends to zero.\\  Thus, the convergence of $R_1$ gives us, as $\varepsilon$ goes to zero,
\begin{equation}
\label{R1 term}
\begin{split}
\left|\int_0^T\int_{\O}R_1\overline{\u}\,dx\,dt\right|&=\left|\int_0^T\int_{\O}\left(\overline{\Dv(\u\otimes\u)}-\Dv(\u\otimes\overline{\u})\right)\overline{u}\,dx\,dt\right|
\\&\leq
\|R_1\|_{L^{\frac{2p}{p+2}}(0,T;L^{\frac{2q}{q+2}}(\O))}\|\overline{\u}\|_{L^p(0,T;L^q(\O))}
\\&\to 0,
\end{split}
\end{equation} for any $p,q\geq 4$.

Letting $\varepsilon$ goes to zero in \eqref{energy equality before limit}, using \eqref{R3 term}, \eqref{R2 term} and \eqref{R1 term}, what we have proved is that in the limit,
\begin{equation*}
\begin{split}&
\int_{\O}|\u|^2\,dx+2\mu\int_0^T\int_{\O}|\nabla\u|^2\,dx\,dt=
\int_{\O}|\u_0|^2\,dx,
\end{split}
\end{equation*}
for any weak solutions with additional condition $\u\in L^p(0,T;L^q(\O))$ with $p\geq 4, q\geq 4.$\\

The final step is to improve the restriction $p,q\geq4.$
Note that, $\u\in L^{\infty}(0,T;L^2(\O))$ and $\u\in L^r(0,T;L^s(\O)),$ thus
\begin{equation*}
\|\u\|_{L^p(0,T;L^q(\O))}\leq C\|\u\|^{\theta}_{L^{\infty}(0,T;L^2(\O))}\|\u\|^{1-\theta}_{L^r(0,T;L^s(\O))},
\end{equation*}
for any $\theta\in(0,1)$ such that
\begin{equation*}
\begin{split}&
\frac{1}{p}=\frac{1-\theta}{r},
\\&\frac{1}{q}=\frac{\theta}{2}+\frac{1-\theta}{s}.
\end{split}
\end{equation*}
This yields
$$\left(\frac{1}{r}+\frac{1}{s}\right)(1-\theta)=\frac{1}{p}+\frac{1}{q}-\frac{\theta}{2}\leq \frac{1}{2}(1-\theta),$$
and hence
$$\frac{1}{r}+\frac{1}{s}\leq \frac{1}{2}$$
with $s\geq 4.$

%\section*{Acknowledgments}

\end{document}